\documentclass[draft, 12pt, a4paper, reqno, oneside]{amsart}
\usepackage{amsmath,amsthm,amssymb,mathrsfs}
\usepackage{color,hyperref,graphicx,comment}
\usepackage[utf8]{inputenc}
\usepackage[all]{xy}
\usepackage{mathtools}
\usepackage{enumitem}
\usepackage{faktor}
\usepackage{tikz-cd}
\usepackage[all]{xy} 
\usepackage{upgreek}
\usepackage{bbold} 
\usepackage{stmaryrd} 
\usepackage{indentfirst}
\usepackage{empheq}
\usepackage{lipsum}
\usepackage{natbib}
\setcitestyle{numbers,square}
\usepackage{geometry}
\numberwithin{equation}{section}

\title{MEAN CURVATURE FLOW SOLITONS FROM SYMMETRY GROUP VIEWPOINT}
\author{XU HAN AND ZHONGHUA HOU}
\address{Dalian University of Technology, Dalian City, China, 116024}
\email{xuhan2016@mail.dlut.edu.cn and zhhou@dlut.edu.cn}

\newcommand{\E}{\ensuremath{\textup{\textbf{E}}}}

\theoremstyle{plain}
\newtheorem{defn}{Definition}
\newtheorem{theorem}{Theorem}
\newtheorem{lemma}{Lemma}

\newtheorem{example}{Example}
\numberwithin{equation}{section}
\numberwithin{defn}{section}
\numberwithin{lemma}{section}
\numberwithin{cor}{section}
\numberwithin{example}{section}

\begin{document}

\begin{abstract}
\noindent
The symmetry group of the mean curvature flow in general ambient Riemannian manifolds is determined, based on which we define generalized solitons to the mean curvature flow. We also provide examples of homothetic solitons in non-Euclidean surfaces and prove that all the affine solutions to the mean curvature flow are self-similar solutions.
\end{abstract}

\maketitle

\section{Introduction}
Since self-similar solutions to the mean curvature flow in Euclidean spaces are significant on several aspects, it is tempting to extend the notion to solitons in general ambient manifolds. Indeed, there have been several successful attempts: Hungerbuehler and Smoczyk \cite{Hun2000} considered the group of isometries of the ambient manifold and constructed examples of rotating solitons; Smoczyk \cite{Smk2001} further studied closed conformal vector fields on the ambient manifold; A. Futaki, K. Hattori and H. Yamamoto \cite{Fut2014} studied mean curvature flow solitons in cone manifolds; J. Alias, J. de Lira, and M. Rigoli \cite{ALR2020} introduced a notion of mean curvature flow soliton using vector fields on the ambient manifold including manifolds of constant sectional curvature, Riemannian products and warped product spaces, and particularly investigated the solitons induced by closed conformal vector fields.

In this paper, we take the symmetry group point of view to define solitons to the mean curvature flow. This viewpoint can be traced back to Richard Hamilton's paper \cite{Ham1995} on page 218, which in the authors' opinion a natural perspective to extend the classical notion, and can be applied to other geometric flows. We can also use this viewpoint to get a refined understanding of certain aspects of the classical self-similar solutions such as Theorem 2 of this paper. Similarly, we can answer a question appears in \cite{Smk2001} on page 176: in an ambient Riemannian manifold, are the conformal solutions in the sense of \cite{Smk2001} able to move along the integral curve of the corresponding conformal vector field? The answer is generally no, unless the conformal vector field is a homothetic one.

We now give the definition of solitons to any geometric flow. The more precise definition of the mean curvature flow solitons is given in Section 5.\par
\begin{defn}
A smooth solution to a geometric flow is called a soliton, if it can be generated by a time-parameter subgroup of the symmetry group acting on the initial hypersurface.
\end{defn}
\noindent

Concerning the symmetry group of the mean curvature flow, Peter Olver \cite{OST1994} computed the symmetry group of the curve shortening flow in the plane. K. Chou and G. Li \cite{Cho2002} obtained the symmetry group of the generalized curve shortening flow in the plane. In this paper, we determine the symmetry group of the mean curvature flow of codimension one in an ambient Riemannian manifold of any dimension.

Let $F :M^{n}\times I \rightarrow (N^{n+1},\bar g)$, where $I$ is an interval, be a smooth family of hypersurface immersions satisfying
\begin{equation}
\label{MCF equation}
\bigg(\frac{\partial F}{\partial t}(x,t)\bigg)^{\bot}=H(x,t)\nu(x,t)\\
\end{equation}
where $H(x,t)$ is the mean curvature with respect to the unit normal $\nu(x,t)$, and $\bot$ denotes the projection along the normal direction.

\begin{theorem}
\label{Th SymGrp}
The infinitesimal symmetries and symmetry transformations of (\ref{MCF equation}) are listed in the following table:

\begin{table}[h]
\centering
\begin{tabular}{|c|c|}
\hline
\textbf{Infinitesimal Symmetry}                         & \textbf{Symmetry Transformation}     \\ \hline
$\partial_t$                                            & Translations in $t$         \\ \hline
${\mathfrak X}(M)$                                      & Diffeomorphisms of $M$      \\ \hline
${\mathcal K}(N)$                                      & Isometries of $N$           \\ \hline
$X+2\lambda t\partial_t$                                & Parabolic Rescalings   \\ \hline
\end{tabular}
\end{table}

\noindent
where ${\mathfrak X}(M)$ and ${\mathcal K}(N)$ denote all the smooth vector fields and Killing vector fields on $M$ respectively, and $X$ is any homothetic vector field on $N$ such that the Lie derivative of the metric $\bar g$ with respect to $X$ satisfies $L_{X}\bar g=2\lambda \bar g$ for a non-zero constant $\lambda$.
\end{theorem}

As an application of Theorem 1, we can solve the problem concerning the classical self-similar solutions in Euclidean spaces. For the curve shortening flow in the plane, Halldorsson \cite{Hal2012} combines the classical self-similar curves in the plane generated by motions such as rotation \cite{Alt1991}, scaling \cite{AbL1986} and translation \cite{Eck2004}. We can define the affine solutions in a similar manner for the  mean curvature flow in Euclidean spaces of any dimension, and prove
\begin{theorem}
All the affine solutions to the mean curvature flow are self-similar solutions.
\end{theorem}
\noindent
\textbf{Remark.} Here we regard the minimal hypersurface, i.e. the stationary solution to the mean curvature flow as a trivial self-similar solution.

Examples of the soliton solutions have been constructed and studied in a wide range of literature, among which we quote \cite{AbL1986}, \cite{ALR2020}, \cite{Eck2004}, \cite{Fut2014}, \cite{Hal2012}, \cite{Hui1990} and \cite{Hun2000}. However, the only homothetic solitons are those in Euclidean spaces which appear in \cite{AbL1986} and \cite{Hui1990}. In this paper, we provide examples of homothetic solitons in non-Euclidean ambient surfaces.

The organization of the paper is the following. In Section 2, we review some basic results on the symmetry group theory to which our main reference is \cite{Olv1993}. In Section 3 and Section 4, we prove Theorem 1 in dimension two and higher dimensions respectively. In Section 5, we formulate the mean curvature flow solitons and construct examples of non-Euclidean homothetic solitons. In Section 6, we prove Theorem 2.

\section{Preliminaries}
Peter Olver \cite{Olv1993} gives a beautiful presentation on the symmetry group theory. For the convenience, we include in this section the necessary concepts and results, and suggest the reader refer to \cite{Olv1993} for the detailed accounts. \par

The basic idea of determining the symmetry group of a PDE system is that we first extend the usual space of variables to the jet space to include partial derivatives, which transforms the original system into a system of algebraic equations, then apply the corresponding results on symmetries of algebraic equations, however, the symmetry of the algebraic system is generally less than that of the original system, so we have to add certain conditions to guarantee the equivalence.

\begin{defn}
Let $M$ be a smooth manifold and $G$ be a Lie group or a local Lie group. Suppose that $\mathcal{U}$ is an open subset of $G\times M$, such that
\begin{equation}
\{e\}\times M\subset \mathcal{U}.
\end{equation}
A local group of transformations acting on $\mathcal{U}$ is defined as a smooth map $\Omega:\mathcal{U}\rightarrow M$ with the following properties:\par
\noindent
(a) If $(h,x)\in\mathcal{U}$, $(g,\Omega(h,x))\in\mathcal{U}$ and $(g\cdot h,x)\in \mathcal{U}$, then
\begin{equation}
\Omega(g,\Omega(h,x))=\Omega(g\cdot h,x).
\end{equation}
(b) For all $x\in \mathcal{U}$,
\begin{equation}
\Omega(e,x)=x.
\end{equation}
(c) If $(g,x)\in \mathcal{U}$, then $(g^{-1}, \Omega(g,x))\in\mathcal{U}$ and
\begin{equation}
\Omega(g^{-1},\Omega(g,x))=x.
\end{equation}
\end{defn}
\noindent
\textbf{Remark.} In this paper, by a one-parameter family of transformations, we mean a one-parameter subgroup of a local group of transformations.

\begin{defn}
Let $G$ be a local group of transformations defined on some open subset $U\subset\{(x,u)\in\mathbb{R}^n\times\mathbb{R}^q\}$, which is the space of independent and dependent variables of a PDE system. We call $G$ is a symmetry group of the system, if $u=f(x)$ is a solution to the system, then whenever $g\cdot f$ is defined for $g\in G$, we have $u=g\cdot f(x)$ is also a solution.
\end{defn}
\noindent
\textbf{Remark.} We only consider the case $q=1$ in the paper, so the variable space becomes $\{(x,u)\in\mathbb{R}^{n+1}\}$ and the independent variable is also denoted by $x=(x^1,...,x^n)$.

\begin{defn}
The $m$-th prolongation of a smooth function $u=f(x)$, which is denoted by $u^{(m)}=pr^{(m)}f(x)$, is an vector-valued function $\{u_{J}\}$ for all multi-indices $J=(j_1,...,j_k)$, and $0\leqslant k\leqslant m$ defined by the following equations
\begin{equation}
u_{J}:=\frac{\partial^{k} f(x)}{\partial x^{j_1}\cdot\cdot\cdot\partial x^{j_k}},
\end{equation}
particularly, we define $u_{J}=u$ when $k=0$.
\end{defn}

\begin{defn}
Given a local group of transformations $G$ acting on an open subset $U$ of the variable space, and for any smooth function $u=f(x)$ such that $u_0=f(x_0)$ and $(\tilde x_0, \tilde u_0)=g\cdot(x_0,u_0)$ is defined for $g\in G$, we define the $m$-th prolongation of $g$ at the point $(x_0,u_0)$ by
\begin{equation}
pr^{(m)}g\cdot(x_0,u_0^{(m)})=(\tilde x_0, \tilde {u}_0^{(m)}),
\end{equation}
where
\begin{equation}
\tilde {u}_0^{(m)}:=pr^{(m)}(g\cdot f)(\tilde x_0).
\end{equation}
\end{defn}

\begin{defn}
Suppose $v$ is a vector field on $U\subset\mathbb{R}^{n+1}$, with the corresponding (local) one-parameter transformations $\Omega_{\varepsilon}$. Then $m$-th prolongation of $v$ at the point $(x,u^{(m)})$ , which is denoted by $pr^{(m)}v(x,u^{(m)})$, is defined by the following equation
\begin{equation}
pr^{(m)}v(x,u^{(m)})=\frac{d}{d\varepsilon}\bigg|_{\varepsilon=0}pr^{(m)}\Omega_{\varepsilon}(x,u^{(m)})
\end{equation}
\end{defn}

\begin{defn}
\label{def max rank}
A system of differential equations
\begin{equation}
\Phi^{\alpha}(x,u^{(m)})=0,\quad \alpha=1,...,l
\end{equation}
is of maximal rank if the Jacobian matrix of $\Phi:=\{\Phi^{\alpha}\}$
\begin{equation}
J_{\Phi}(x,u^{(m)})=\{(\frac{\partial \Phi^{\alpha}}{\partial x^{i}}),(\frac{\partial \Phi^{\alpha}}{\partial u_{J}})\},
\end{equation}
is of rank $l$ on $\mathcal{S}:=\{(x,u^{(m)}):\Phi(x,u^{(m)})=0\}$.
\end{defn}

\begin{defn}
\label{def loca solv}
A system of differential equations $\Phi(x,u^{(m)})=0$ is said to be locally solvable at a point $(x_{0}, u_{0}^{(m)})\in \mathcal{S}$, if there is a solution $u=f(x)$ of the system defined in a neighbourhood of $x_{0}$, such that $u_{0}^{(m)}=pr^{(m)}f(x_{0})$

\end{defn}

\begin{defn}
A system of differential equations is called non-degenerate, if it is both locally solvable and of maximal rank at every point $(x_{0}, u_{0}^{(m)})\in \mathcal{S}$.

\end{defn}

We have the following necessary and sufficient condition for a group to be the symmetry group.

\begin{theorem}
\label{equi cond}
Let
\begin{equation}
\Phi^{\alpha}(x,u^{(m)})=0,\quad \alpha=1,...,l
\end{equation}
be a non-degenerate system of differential equations defined on $U\subset \mathbb{R}^{n+1}$. If $G$ is a connected local group of transformations acting on $U$, then $G$ is a symmetry group of the system if and only if for every infinitesimal generator $v$ of $G$, the following are satisfied on $\mathcal{S}$,
\begin{equation}
\label{theorem 2 cond}
pr^{(m)}v[\Phi^{\alpha}(x,u^{(m)})]=0,
\end{equation}
for $\alpha=1,...,l$.
\end{theorem}

We also need a prolongation formula.
\begin{defn}
The $i$-th total derivative of a given function $P(x,u^{(m)})$ is defined by the following equation
\begin{equation}
D_iP=\frac{\partial P}{\partial x^i}+\sum_J u_{J,i}\frac{\partial P}{\partial u_J},
\end{equation}
where $J=(j_1,...,j_k)$, and
\begin{equation}
u_{J,i}=\frac{\partial u_J}{\partial x^i}=\frac{\partial^{k+1} u}{\partial x^i \partial x^{j_1}\cdot\cdot\cdot\partial x^{j_k}}.
\end{equation}
The $J$-th total derivative is defined by
\begin{equation}
D_J=D_{j_1}D_{j_2}\cdot\cdot\cdot D_{j_k}.
\end{equation}
\end{defn}

\begin{theorem}
\label{formula}
Let
\begin{equation}
v=\sum_{i=1}^{n}\xi^{i}(x,u)\frac{\partial}{\partial x^{i}}+\eta(x,u)\frac{\partial}{\partial u}
\end{equation}
be a vector field defined on the variable space. Suppose that the $m$-th prolongation of $v$ is
\begin{equation}
pr^{(m)}v=v+\sum_{J}\phi^{J}(x,u^{(m)})\frac{\partial}{\partial u_{J}}.
\end{equation}
Then the coefficient functions are given by
\begin{equation}
\phi^{J}(x,u^{(m)})=D_{J}\bigg(\eta-\sum_{i=1}^{n}\xi^{i}u_{i}\bigg)+\sum_{i=1}^{n}\xi^{i}u_{J,i}.
\end{equation}
\end{theorem}

\section{Proof of Theorem 1 in dimension two}
In this section, we prove Theorem \ref{Th SymGrp} when the ambient manifolds are surfaces. And in Section 4, we apply the same procedures to all dimensions. The result on surfaces is a little stronger than that of higher dimensions due to technical reasons that will be seen along the proof. \par

The outline of the proof is given as follows: firstly, we combine the mean curvature flow equation and the metric evolution equation to form a system, and use Theorem \ref{formula} and the sufficient part of Theorem \ref{equi cond} to calculate the determining equations; secondly, in order to apply the necessary part of Theorem \ref{equi cond}, we have to check the non-degeneracy condition of the system. Since the symmetry group of the system is exactly that of the mean curvature flow, we complete the proof.

Let $N$ be a Riemannian surface with metric $\bar g$ and Levi-Civita connection $D$. Suppose that a smooth family of immersed curves $F: M\times I \rightarrow (N,\bar g)$ is a solution to the mean curvature flow such that
\begin{equation}
\label{CSF equation}
\bigg(\frac{\partial F}{\partial t}(x,t)\bigg)^{\bot}=k(x,t)\nu(x,t),
\end{equation}
where $k(x,t)$ is the geodesic curvature of the curves in $N$ with respect to the unit normal $\nu(x,t)$. The induced metric of $M$ is denoted by $g$ and the corresponding Levi-Civita connection is denoted by $\nabla$. We also define notations $F_t(x):=F(x,t)$ and $M_t=F_t(M)$.\par

\begin{lemma}(Theorem 1.3.13 of \cite{Ger2006})
\label{Gauss coordinate}
Let $(N^{n+1},\bar g)$ be a Riemannian manifold and $M\subset N$ be a connected smooth hypersurface. Then for any point $p\in M$, there exists a neighbourhood $U\subset N$ and a coordinate system $\{x^\alpha\}$, $0\leqslant\alpha\leqslant n$ such that
\begin{equation}
M\cap U=\{x^0=0\}.
\end{equation}
and
\begin{equation}
\bar g|_{U}=(dx^{0})^2+\sum_{i,j=1}^{n}\sigma_{ij}(dx^i)(dx^j).
\end{equation}
The local chart $(U,\{x^\alpha\})$ is usually called a normal Gaussian coordinate system.
\end{lemma}
\noindent
By Lemma \ref{Gauss coordinate}, we can choose a normal Gaussian coordinate system $\{x,y\}$ of a neighbourhood $U\subset N$ containing the initial curve locally, such that
\begin{equation}
\bar g=A(x,y)dx^2+dy^2.
\end{equation}
We can further choose $x$ to be the arc-length parameter of the initial curve such that $A(x,0)=1$. By the following lemma, we can represent $M_t$ locally as graphs. Since the computation of the symmetries is in local, no generality is lost under the assumptions.

\begin{lemma}
\label{graph lemma}
Let $F(x,t)$ be a solution to the mean curvature flow. During some time interval $I'$ short enough, for each $t\in I'$, $M_t$ is locally the graph of a function $u(x,t)$.
\end{lemma}
\noindent
Proof. For any $t_0\in I$, by Lemma \ref{Gauss coordinate} we can choose a normal Gaussian coordinates $\{x^{\alpha}\}$ of some neighbourhood $U\subset N$ such $M_{t_0}\cap U=\{x^0=0\}$ and $x=(x^1,...,x^n)$ is the local coordinates of $M_{t_0}$. Since $\nabla x_0=0$ on $M_{t_0}$, there exists a time interval $I'$ short enough such that $\nabla x_0(x,t)$ is bounded on $M_t\cap U$ for any $t\in I'$. Thus $M_t\cap U$ is the graph of a function $x^0=u(x,t)$.\par
\hfill\qedsymbol

The tangent vector field to the curve is
\begin{equation}
F_*(\frac{\partial}{\partial x})=\frac{\partial}{\partial x} +u_x\frac{\partial}{\partial y},
\end{equation}
and the induced metric $g$ on $M$ is then
\begin{equation}
g=L^2dx^2,
\end{equation}
where $L^2:=A+u_x^2$. We choose the unit normal to be
\begin{equation}
\nu=\frac{1}{\sqrt{A}L}(-u_{x}\partial x +A\partial y),
\end{equation}
and the corresponding geodesic curvature is given by
\begin{equation}
\label{2dim curv}
k=\frac{\sqrt{A}}{L^{3}}\Big(u_{xx}-\frac{A_y}{A}u_{x}^2-\frac{A_x}{2A}u_{x}-\frac{1}{2}A_y\Big),
\end{equation}
where $A_x$ and $A_y$ are partial derivatives of $A$. The mean curvature flow equation under the graph representation can be written as
\begin{equation}
\bar g( \frac{\partial F}{\partial t}, \nu )=\frac{\sqrt{A}}{L}u_t=k,
\end{equation}
which is equivalent to
\begin{equation}
\label{2-dim graph equ}
u_t=\frac{1}{L^2}\Big(u_{xx}-\frac{A_y}{A}u_{x}^2-\frac{A_x}{2A}u_{x}-\frac{1}{2}A_y\Big),
\end{equation}
and we set
\begin{equation}
\Phi^1:=L^2u_t-u_{xx}+\frac{A_y}{A}u_{x}^2+\frac{A_x}{2A}u_{x}+\frac{1}{2}A_y.
\end{equation}\par

The evolution equation of the induced metric under graph representation is given by (cf. \cite{GaH1986})
\begin{equation}
\frac{\partial L^2}{\partial t}=-2L^2k^2,
\end{equation}
that is
\begin{equation}
\label{metric evol 2d}
u_xu_{xt}=-L^2k^2,
\end{equation}
and we set
\begin{equation}
\Phi^2:=u_xu_{xt}+L^2k^2.
\end{equation}
Combining the two equations (\ref{2-dim graph equ}) and (\ref{metric evol 2d}) into a system
\begin{equation}
\label{couple equ 2d}
\begin{cases}
\Phi^1=0\\
\Phi^2=0.
\end{cases}
\end{equation}
Generally, we need to find the infinitesimal symmetries of (\ref{couple equ 2d}) in the form of
\begin{equation}
v:=\tau(t,x,u)\frac{\partial}{\partial t}+\xi(t,x,u)\frac{\partial}{\partial x}+\eta(t,x,u)\frac{\partial}{\partial u},
\end{equation}
but for the geometric flow, which is invariant under all diffeomorphisms of $M$, we can choose an arbitrary fixed coordinates of $M$ such that $\xi=\xi(x,u)$ and $\tau=\tau(t)$. In fact, given any fixed coordinates $(x,u)$ of $U\subset N$ containing $M$ locally, and suppose $g_\varepsilon$ preserves the solution
$u(x,t)$ as graphs for some short time interval, then the vector field of
\begin{equation}
g_\varepsilon(t,x,u(x,t)):=(t_\varepsilon,x_\varepsilon,u_\varepsilon(x_\varepsilon,t_\varepsilon)),
\end{equation}
is
\begin{equation}
\frac{d}{d\varepsilon}\bigg|_{\varepsilon=0}g_\varepsilon(t,x,u(x,t))
=\frac{d t_\varepsilon}{d\varepsilon}\bigg|_{\varepsilon=0}\frac{\partial}{\partial t}+\frac{dx_\varepsilon}{d\varepsilon}\bigg|_{\varepsilon=0}\frac{\partial}{\partial x}+\frac{du_\varepsilon}{d\varepsilon}\bigg|_{\varepsilon=0}\frac{\partial}{\partial u}.
\end{equation}
Since $x$ and $t$ are independent variables, we prove that $\xi=\xi(x,u)$ and $\tau=\tau(t)$. In dimension two this can also be easily derived from the determining equations, which can be seen along the proof of Theorem 1.

Proof of Theorem 1:

\textbf{Step 1.} Applying Theorem \ref{formula} to $v$, we obtain the prolongation of the vector $v$
\begin{equation}
\label{prolong v}
pr^{(2)}v=v+\phi^x\frac{\partial}{\partial u_x}+\phi^t\frac{\partial}{\partial u_t}+\phi^{xx}\frac{\partial}{\partial u_{xx}}+\phi^{xt}\frac{\partial}{\partial u_{xt}}+\phi^{tt}\frac{\partial}{\partial u_{tt}},
\end{equation}
where the coefficients are given by the the formula
\begin{equation}
\begin{aligned}
\phi^x=&\eta_x+(\eta_u-\xi_x)u_x-\tau_xu_t-\xi_uu_x^2-\tau_uu_xu_t\\
\phi^t=&\eta_t-\xi_tu_x+(\eta_u-\tau_t)u_t-\xi_uu_xu_t-\tau_uu_t^2\\
\phi^{xx}=&\eta_{xx}+2\eta_{xu}u_x-\xi_{xx}u_x-\tau_{xx}u_t+\eta_{uu}u_x^2-2\xi_{xu}u_x^2-2\tau_{xu}u_xu_t-\xi_{uu}u_x^3\\
          &-\tau_{uu}u_x^2u_t+\eta_uu_{xx}-2\xi_xu_{xx}-2\tau_xu_{xt}-3\xi_uu_xu_{xx}-\tau_uu_tu_{xx}-2\tau_uu_xu_{xt}\\
\phi^{xt}=&\eta_{tx}+\eta_{tu}u_x+\eta_{ux}u_t-\xi_{tx}u_x-\tau_{tx}u_t+\eta_{uu}u_xu_t-\xi_{tu}u_x^2-\xi_{ux}u_tu_x\\
          &-\tau_{tu}u_xu_t+\eta_uu_{tx}-\xi_tu_{xx}-\xi_xu_{xt}-\tau_tu_{tx}-\xi_{uu}u_tu_x^2-\tau_{ux}u_t^2\\
          &-2\xi_uu_{tx}u_x-\xi_uu_tu_{xx}-\tau_{uu}u_xu_t^2-2\tau_uu_tu_{tx}-\tau_xu_{tt}-\tau_uu_xu_{tt}.
\end{aligned}
\end{equation}
By the sufficient condition of Theorem \ref{equi cond}, we apply $pr^{(2)}v$ to (\ref{couple equ 2d}) to obtain
\begin{equation}
\label{prv=0 2d}
\begin{cases}
pr^{(2)}v[\Phi^1(x,u^{(2)})]=0\\
pr^{(2)}v[\Phi^2(x,u^{(2)})]=0.
\end{cases}
\end{equation}
which is satisfied whenever $\Phi^1=0$ and $\Phi^2=0$. This can be regarded as a system of algebraic equations where the dependence of monomials is further restricted to the condition $\Phi^1=0$ and $\Phi^2=0$. \par

An observation is that the second equation of (\ref{prv=0 2d}) is immaterial for calculating the infinitesimal symmetries. By straightforward calculation, we obtain that $pr^{(2)}v[\Phi^1(x,u^{(2)})]=0$ is equivalent to
\begin{equation}
\label{prv(phi1)=0}
\begin{aligned}
0=&u_{xt}(2u_x\tau_u+2\tau_x)-u_{t}^2(A\tau_u+3u_x^2\tau_u+2u_x\tau_x)\\
  &+u_{xx}(3u_x\xi_u+2\xi_x-\eta_u)+u_tF_1(u_x)+F_2(u_x).
\end{aligned}
\end{equation}
where $F_1$ and $F_2$ are both polynomials with respect to $u_x$.
It can be seen that $\Phi^1=0$ and $\Phi^2=0$ only provide extra dependence of $u_t$ and $u_{xt}$ as follows
\begin{equation}
\label{depend relation}
\begin{cases}
u_t=\frac{L}{\sqrt{A}}k\\
u_xu_{xt}=-L^2k^2.
\end{cases}
\end{equation}
Inserting (\ref{depend relation}) into (\ref{prv(phi1)=0}), we find that among all independent monomials containing $u_{xx}^2$ there are two terms
\begin{equation}
-3\tau_u\frac{A}{L^4}u_{xx}^2 \quad and \quad -2\tau_x\frac{u_xu_{xx}^2}{L^4}
\end{equation}
This implies $\tau_u=0$ and $\tau_x=0$ due to $A>0$. Thus (\ref{prv(phi1)=0}) contains no terms of $u_{xt}$, which implies that $\Phi^2=0$ provides no further restrictions on (\ref{prv(phi1)=0}).\par

\textbf{Step 2.} We now focus on $pr^{(2)}v[\Phi^1(x,u^{(2)})]=0$ to calculate the infinitesimal symmetries. By expanding (\ref{prv(phi1)=0}) and considering $\tau=\tau(t)$, we get the following PDE system of determining equations by some simple reductions
\begin{subequations}
\label{determing equation 2d}
\begin{empheq}[left=\empheqlbrace]{align}
&2A\xi_x+A_x\xi+A_u\eta =A\tau_t\\
&A\xi_u+\eta_x=0\\
&2\eta_u=\tau_t\\
&(A\xi_u)_u=A\xi_t \\
&\big(\frac{A_u}{A}\big)_u\eta+\big(\frac{A_u}{A}\big)_x\xi+\frac{A_x}{A}\xi_u+2\xi_{xu}+\eta_t+\frac{A_u}{A}\eta_u=0\\
&\big(\frac{A_x}{A}\big)_u\eta+\big(\frac{A_x}{A}\big)_x\xi-2A\xi_t+\frac{A_u}{A}\eta_x+\frac{A_x}{A}\xi_x+2\xi_{xx}=0\\
&A_{uu}\eta+A_{ux}\xi+2A_u\xi_x+2A\eta_t-A_u\eta_u+\frac{A_x}{A}\eta_x-2\eta_{xx}=0\\
&\tau=\tau(t).
\end{empheq}
\end{subequations}
We claim that $\eta$ and $\xi$ are independent on $t$. By (\ref{determing equation 2d}c) and (\ref{determing equation 2d}h), we have $\eta_u$ is independent on $u$ and $x$ , thus $\eta=\frac{\tau_t}{2}u+\beta(x,t)$ for some function $\beta(x,t)$. And by (\ref{determing equation 2d}b), we have $A\xi_u$ is independent on $u$, thus $A\xi_u=\alpha(x,t)$ for some function $\alpha(x,t)$, then by (\ref{determing equation 2d}d), we have $\xi=\xi(x,u)$. Differentiating (\ref{determing equation 2d}a) with respect to $u$ and using (\ref{determing equation 2d}e), we can finally obtain $\eta=\eta(x,u)$.\par

Since the left hand side of (\ref{determing equation 2d}c) does not depend on $t$ and its right hand side depends only on $t$, we know that $\tau_t$ must be a constant, so $\tau=2\lambda t+c_1$ for some constants $\lambda$ and $c_1$. Therefore, (\ref{determing equation 2d}a),(\ref{determing equation 2d}b) and (\ref{determing equation 2d}c) can be further reduced to
\begin{equation}
\label{conf kill equa}
\begin{cases}
2A\xi_x+A_x\xi+A_u\eta =2A\lambda \\
A\xi_u+\eta_x=0 \\
\eta_u=\lambda.
\end{cases}
\end{equation}
In fact, the above equations are exactly the characteristic equations of homothetic vector fields on $N$. Suppose that
\begin{equation}
X:=\xi(x,u)\frac{\partial}{\partial x}+\eta(x,u)\frac{\partial}{\partial u}
\end{equation}
is a conformal vector field on $N$ such that $L_X\bar g=2\lambda \bar g$, then we can check that $\xi$ and $\eta$ satisfy (\ref{conf kill equa}) by the formula $\bar g(D_W X, Y )+\bar g(D_Y X, W)=L_X\bar g(W,Y)$ for arbitrary vector fields $W$ and $Y$ on $N$.\par

Since $\xi(x,u)$ and $\eta(x,u))$ satisfying (\ref{conf kill equa}) and $\tau=2\lambda t+c_1$, when $\lambda=0$, $c_1=1$ and $\xi(x,u)=\eta(x,u)=0$, the infinitesimal symmetry is $\partial_t$ corresponding to a translation in $t$. When $\lambda=c_1=0$, we obtain the Killing vector fields corresponding to the isometries of $N$. And when $\lambda\neq 0$ and $c_1=0$, the infinitesimal symmetry is $X+2\lambda t\partial_t$ corresponding to the parabolic rescalings. It is easy to check that all these transformations indeed preserve the mean curvature flow except maybe the parabolic rescalings. Let $\Omega_{\varepsilon}$ be the corresponding transformation acting on $(t,x,u)$ such that

\begin{equation}
\frac{d \Omega_{\varepsilon}(t,x,u)}{d \varepsilon}\bigg|_{\varepsilon=0}=X+2\lambda t\partial_t.
\end{equation}
If we assume

\begin{equation}
\begin{cases}
t(\varepsilon):=\varphi_{\varepsilon}(t)  \\
(x(\varepsilon), u(\varepsilon)):=\omega_{\varepsilon}(x,u)
\end{cases}
\end{equation}
then
\begin{equation}
\frac{dt(\varepsilon)}{d\varepsilon}=2\lambda t(\varepsilon).
\end{equation}
Noting that $\varphi_{0}(t)=t$, we obtain
\begin{equation}
\varphi_{\varepsilon}(t)=e^{2\lambda\varepsilon}t.
\end{equation}
Thus
\begin{equation}
\label{action 2d}
\Omega_{\varepsilon}(t,x,u)=(e^{2\lambda\varepsilon}t,\omega_{\varepsilon}(x,u)),
\end{equation}
where
\begin{equation}
\frac{d \omega_{\varepsilon}(x,u)}{d \varepsilon}\bigg|_{\varepsilon=0}=X.
\end{equation}
Since $X$ is a homothetic vector field, we see that $\omega_{\varepsilon}$ acting on $N$ is a a smooth family of homothetic transformations. We can calculate the conformal factor as follows.
\begin{lemma}
\label{conformal factor}
Suppose $\omega_\varepsilon^*\bar g=c^2(\varepsilon)\bar g$ and $c(\varepsilon)>0$, then $c(\varepsilon)=e^{\lambda \varepsilon}$.
\end{lemma}
\noindent
Proof. This is directly from the definition of Lie derivative. Since
\begin{equation}
(L_X)^k\bar g=\lim\limits_{\varepsilon\rightarrow 0}\frac{(\omega_\varepsilon^*-Id)^k\bar g}{\varepsilon^k}=(2\lambda)^k\bar g,
\end{equation}
and by induction on $k$ we have
\begin{equation}
\lim\limits_{\varepsilon\rightarrow 0}\frac{(\omega_\varepsilon^*-Id)^k\bar g}{\varepsilon^k}=\frac{d^k(\omega_\varepsilon^*\bar g)}{d\varepsilon^k}\bigg|_{\varepsilon=0}.
\end{equation}
We complete the proof using Taylor's expansion of $\omega_\varepsilon^*\bar g$
\begin{equation}
\omega_\varepsilon^*\bar g=\sum_{k=0}^{\infty}\frac{\varepsilon^k}{k!}\frac{d^k(\omega_\varepsilon^*\bar g)}{d\varepsilon^k}\bigg|_{\varepsilon=0}\bar g=e^{2\lambda \varepsilon}\bar g.
\end{equation}
$\hfill\qedsymbol$\par
Since the homothetic transformations $\omega_{\varepsilon}$ of $N$ can be regarded as equipping $N$ with a smooth family of metric $\omega_{\varepsilon}^*\bar g$, we obtain by Lemma \ref{conformal factor}
\begin{equation}
\nu|_{\omega_{\varepsilon}^*\bar g}(x,u)=e^{-\lambda\varepsilon}\nu|_{\bar g}(x,u),
\end{equation}
\begin{equation}
H|_{\omega_{\varepsilon}^*\bar g}(x,u)=e^{-\lambda\varepsilon}H|_{\bar g}(x,u),
\end{equation}
that is
\begin{equation}
(H\nu)|_{\omega_{\varepsilon}^*\bar g}=e^{-2\lambda\varepsilon}(H\nu)|_{\bar g}.
\end{equation}
Therefore, considering (\ref{action 2d}), $\Omega_{\varepsilon}$ indeed preserve the mean curvature flow. In summary, besides the diffeomorphisms we have four types of infinitesimal symmetries which span the Lie algebra of the symmetry group.

\textbf{Step 3.} In order to apply Theorem \ref{equi cond} to complete the proof, it remains to check the non-degeneracy condition. Since
\begin{equation}
\label{mcf 2dim}
\begin{cases}
\Phi^1=L^2\big(u_t-\frac{L}{\sqrt{A}}k\big)\\
\Phi^2=u_xu_{xt}+L^2k^2,
\end{cases}
\end{equation}
we have
\begin{equation}
\frac{\partial\Phi^1}{\partial u_t}=L^2.
\end{equation}
and
\begin{equation}
J_{\Phi^2}(x,t,u^{(2)})\equiv0
\end{equation}
on $\{(x,t,u^{(2)}):\Phi^1=0, \Phi^2=0\}$ if and only if $u_x\equiv0$, i.e. $k\equiv0$. By Definition \ref{def max rank}, the system is of maximal rank whenever the solution is not a geodesic.\par
By Definition \ref{def loca solv}, in order to check the local solvability condition, we have to find a solution $u=f(x,t)$ of the flow such that $u_0^{(2)}=pr^{(2)}f(x_0,t_0)$, for any given data \begin{equation}
(x_0, t_0, u_0, (u_t)_0, (u_x)_0, (u_{xx})_{0}, (u_{xt})_{0}, (u_{tt})_{0}),
\end{equation} satisfying
\begin{equation}
\label{initial restriction}
\begin{cases}
(u_t)_0=\frac{L_0}{\sqrt{A_0}}k_0\\
(u_x)_0(u_{xt})_{0}=-L_0^2k_0^2,
\end{cases}
\end{equation}
where $(u_t)_0:=u_t(x_0, t_0)$, $A_0:=A(x_0, u_0)$, $L_0:=L(x_0, u_0)$ and $k_0:=k(x_0, u_0)$. Since the flow is invariant under translations in $t$, we can assume $t=0$. And it suffices to construct a initial curve $f(x, 0)=f_0(x)$ such that $u_0^{(2)}=pr^{(2)}f_0(x_0)$. Then we can involve the initial curve to construct a solution $f(x,t)$ to the mean curvature flow due to short-time existence of the solution. \par

We assume that
\begin{equation}
f_0(x)=u_0+(u_x)_0(x-x_0)+\frac{1}{2!}(u_{xx})_{0}(x-x_0)^2+\frac{1}{4!}C(x-x_0)^4,
\end{equation}
where $C$ is a constant to be determined, and $u=f(x,t)$ satisfies

\begin{equation}
\label{f(x,t)sys}
\begin{cases}
f_t-\frac{L}{\sqrt{A}}k=0\\
f_xf_{xt}+L^2k^2=0\\
f(x,0)=f_0(x).
\end{cases}
\end{equation}
It can be easily seen from the definition of $f_0(x)$ that
\begin{equation}
\begin{gathered}
f(x_0,0)=u_0,
\quad
f_x(x_0,0)=(u_x)_0,
\quad
f_{xx}(x_0,0)=(u_{xx})_0.
\end{gathered}
\end{equation}
Since $f_t-\frac{L}{\sqrt{A}}k=0$, we see that $f_t(x_0,0)$ is completely determined by $x_0$, $u_0$, $(u_x)_0$ and $(u_{xx})_{0}$, and by the first equation of (\ref{initial restriction}), it is easy to check that
\begin{equation}
f_t(x_0,0)=(u_t)_0.
\end{equation}

For the initial value $(u_{xt})_{0}$, we have to consider the evolution of the metric
\begin{equation}
f_xf_{xt}+L^2k^2=0.
\end{equation}
\noindent
Similarly we see that $f_{xt}(x_0,0)$ is also determined by $x_0$, $u_0$, $(u_x)_0$ and $(u_{xx})_{0}$ and by the second equation of (\ref{initial restriction}), we have
\begin{equation}
f_{xt}(x_0,0)=(u_{xt})_0.
\end{equation}

For the initial value $(u_{tt})_{0}$, we differentiate $u_t-\frac{L}{\sqrt{A}}k$ with respect to $t$
\begin{equation}
\label{expr utt}
u_{tt}=\Big(\frac{L}{\sqrt{A}}\Big)_tk+\frac{L}{\sqrt{A}}k_t.
\end{equation}
and use the evolution equation of the curvature (cf. Lemma 10.7 of \cite{Zhu2002})
\begin{equation}
\frac{\partial k}{\partial t}=\frac{1}{L^2}\frac{\partial^2 k}{\partial x^2}+(k^2+R)k,
\end{equation}
where $R$ is the Gaussian curvature of $N$. By the expression of $k$, i.e. (\ref{2dim curv}), we obtain
\begin{equation}
\label{expr kt}
\frac{\partial k}{\partial t}=\frac{\sqrt A}{L^5}u_{xxxx}+\alpha(u_{xxx}, u_{xx}, u_x, x, u)+\beta(u_{xx}, u_x, x, u),
\end{equation}
for some function $\alpha$ and $\beta$. Inserting (\ref{expr kt}) into (\ref{expr utt}), and noting that
\begin{equation}
\Big(\frac{L}{\sqrt{A}}\Big)_tk=\gamma(u_{xx}, u_x, x, u),
\end{equation}
for some function $\gamma$, we finally obtain
\begin{equation}
\begin{aligned}
f_{tt}(x_0,0)=&\gamma((u_{xx})_{0}, (u_x)_0, x_0, u_0)\\
              &+\frac{C}{L_0^4}+\frac{L_0}{\sqrt{A_0}}\Big(\alpha(0,(u_{xx})_{0}, (u_x)_0, x_0, u_0)+\beta((u_{xx})_{0}, (u_x)_0, x_0, u_0)\Big)
\end{aligned}
\end{equation}
and we can solve $C$ such that $f_{tt}(x_0,0)=(u_{tt})_{0}$ due to $L_0\neq0$. \par

$\hfill\qedsymbol$

\section{Proof of Theorem 1 in any dimensions}
In this section, we set the index notations: $0\leqslant\alpha, \beta, \gamma,...\leqslant n$ and $1\leqslant i, j, k,...\leqslant n$ and use Einstein summation convention unless otherwise stated.\par

Suppose that a smooth family of hypersurface immersions $F :M^{n}\times I \rightarrow (N^{n+1},\bar g)$ satisfies the mean curvature flow equation
\begin{equation}
\label{MCF equ}
\bigg(\frac{\partial F}{\partial t}(x,t)\bigg)^{\bot}=H(x,t)\nu(x,t),
\end{equation}

By the same arguments as those in Section 3, we can choose a normal Gaussian coordinate system $\{x^\alpha\}$ of $U\subset N$ containing $F_0(M)$ locally such that

\begin{equation}
\bar g=(dx^0)^2+\sigma_{ij}(x^0,x)dx^idx^j,
\end{equation}
and the family of hypersurfaces can be represented as graphs $x^0=u(x^1,...,x^n,t)$ during some short time interval. The induced metric on $M$ can be written as
\begin{equation}
g_{ij}=u_iu_j+\sigma_{ij},
\end{equation}
where $u_i$ is the partial derivative of $u$ with respect to $x^i$.
And its inverse is
\begin{equation}
g^{ij}=\sigma^{ij}-\frac{u^iu^j}{L^2},
\end{equation}
where $\sigma^{ij}$ is the inverse of $\sigma_{ij}$, $u^i=\sigma^{ij}u_j$ and $L^2=1+\sigma^{ij}u_iu_j$.
The unit normal of the graph is
\begin{equation}
\nu=(1,-u^1,...,-u^n)\frac{1}{L}.
\end{equation}
By the Gauss formula
\begin{equation}
D_{\frac{\partial}{\partial x^i}}\frac{\partial}{\partial x^j}=\nabla_{\frac{\partial}{\partial x^i}}\frac{\partial}{\partial x^j}+h(\frac{\partial}{\partial x^i},\frac{\partial}{\partial x^j})\nu,
\end{equation}
the second fundamental form $h_{ij}=h(\frac{\partial}{\partial x^i},\frac{\partial}{\partial x^j})$ of $M$ satisfies
\begin{equation}
\frac{\partial^2 F^\alpha}{\partial x^i \partial x^j}+\Gamma_{\beta\gamma}^\alpha\frac{\partial F^\beta}{\partial x^i}\frac{\partial F^\gamma}{\partial x^j}-\Gamma_{ij}^{k}\frac{\partial F^\alpha}{\partial x^k}=h_{ij}\nu^\alpha.
\end{equation}
Then
\begin{equation}
h_{ij}=\frac{1}{L}\{u_{ij}-u_l(\Gamma_{0j}^{l}u_i+\Gamma_{i0}^{l}u_j+\Gamma_{ij}^{l})+\Gamma_{ij}^{0}\}.
\end{equation}
The mean curvature follows by taking the trace
\begin{equation}
H=\frac{1}{L}g^{ij}\{u_{ij}-u_l(\Gamma_{0j}^{l}u_i+\Gamma_{i0}^{l}u_j+\Gamma_{ij}^{l})+\Gamma_{ij}^{0}\}.
\end{equation}
Since
\begin{equation}
\bar g( \frac{\partial F}{\partial t}, \nu)=\frac{u_t}{L}=H,
\end{equation}
we obtain
\begin{equation}
\Psi^1:=u_tL^2-(L^2\sigma^{ij}-u^iu^j)(u_{ij}-\Gamma_{0j}^{l}u_iu_l-\Gamma_{0i}^{l}u_ju_l-\Gamma_{ij}^{l}u_l+\Gamma_{ij}^{0})=0.
\end{equation}

Under the graph representation, the evolution equation of the induced metric
\begin{equation}
\frac{\partial g_{kl}}{\partial t}=-2Hh_{kl}
\end{equation}
becomes
\begin{equation}
\begin{aligned}
\Psi^{kl}:=&L^2(u_{kt}u_l+u_{lt}u_k)+2g^{ij}(u_{ij}-\Gamma_{0j}^{s}u_iu_s-\Gamma_{0i}^{s}u_ju_s-\Gamma_{ij}^{s}u_s+\Gamma_{ij}^{0})(u_{kl}\\
                         &-\Gamma_{0l}^{s}u_ku_s-\Gamma_{k0}^{s}u_lu_s-\Gamma_{kl}^{s}u_s+\Gamma_{kl}^{0})=0.
\end{aligned}
\end{equation}

Proof of Theorem 1. The principles of calculation are exactly the same as those in the two dimensional case. And by the arguments in Section 3, it suffices to find the infinitesimal symmetries of the form
\begin{equation}
v:=\tau(t)\frac{\partial}{\partial t}+\xi^i(x,u)\frac{\partial}{\partial x^i}+\eta(t,x,u)\frac{\partial}{\partial u}.
\end{equation}

\textbf{Step 1.} By Theorem \ref{formula}, the prolongation of $X$ is
\begin{equation}
pr^{(2)}v=v+\eta^i\frac{\partial}{\partial u_i}+\eta^t\frac{\partial}{\partial u_t}+\eta^{ij}\frac{\partial}{\partial u_{ij}}+\eta^{it}\frac{\partial}{\partial u_{it}}+\eta^{tt}\frac{\partial}{\partial u_{tt}},
\end{equation}
where the coefficients are computed similarly by the formula in Theorem \ref{formula}. By Theorem \ref{equi cond}, we apply $pr^{(2)}v$ to $\Psi^1$ and $\Psi^{kl}$
\begin{equation}
\begin{cases}
pr^{(2)}v[\Psi^1(x,t,u^{(2)})]=0\\
pr^{(2)}v[\Psi^{kl}(x,t,u^{(2)})]=0.
\end{cases}
\end{equation}
whenever $\Psi^1=0$ and $\Psi^{kl}=0$. Inserting $\Psi^1=0$ into $pr^{(2)}v[\Psi^1(x,t,u^{(2)})]=0$, all the terms containing $u_{jt}$ are $L^2g^{ij}u_{tj}(\tau_i+u_i\tau_u)$, which equal zero due to $\tau=\tau(t)$. Thus $\Psi^{kl}=0$ provides no further dependence.

\textbf{Step 2.} To calculate the determining equations, we only have to consider the first equation $pr^{(2)}v[\Psi^1(x,t,u^{(2)})]=0$. We need not find all the independent monomials and their coefficients. Actually in order to prove Theorem 1, it suffices to focus on the terms containing $u_t$ and the second order derivatives of $u$. This is the main trick that makes the tedious calculation in higher dimensions relatively easy to handle. \par

We denote $L^2g^{ij}$ by $G^{ij}$, then the expansion of $pr^{(2)}v[\Psi^1(x,t,u^{(2)})]=0$ can be written as
\begin{equation}
\begin{aligned}
0=&2G^{kl}u_{ql}(\xi^q_k+u_k\xi^q_u)+u_{kl}G^{kl}(u_q\xi^q_u-\eta_u)   \\
&+u_{kl}\big[-\xi^qG^{kl}_q-\eta G^{kl}_u+u_p\xi^p_qG^{kl}_{u_q}+u_pu_q\xi^p_uG^{kl}_{u_q}-\eta_qG^{kl}_{u_q}-u_q\eta_uG^{kl}_{u_q}\big]   \\
&+u_t\big[\eta(L^2)_u-u_p(L^2)_{u_q}\xi^p_q-u_pu_q(L^2)_{u_q}\xi^p_u+(L^2)_p\xi^p-L^2u_p\xi^p_u-L^2\tau_t   \\
&+(L^2)_{u_p}\eta_p+L^2\eta_u+u_p(L^2)_{u_p}\eta_u\big]+F_0(u_1,...,u_n),
\end{aligned}
\end{equation}
where $(L^2)_{u_q}$ and $G^{kl}_{u_q}$ is partial derivatives with respect to $u_q$, and $F_0(u_1,...,u_n)$ is polynomial with respect to $u_1,...,u_n$. \par

Inserting the dependence relation $\Psi^1=0$ and collect again the second order derivatives of $u$, we find that the coefficient of $u_{kl}$ for each $(k,l)$ is

\begin{equation}
\label{0 order term}
-\xi^s\sigma_{s}^{kl}-\eta\sigma_{u}^{kl}+\xi_{s}^{k}\sigma^{sl}+\xi_{s}^{l}\sigma^{ks}-\tau_t\sigma^{kl}=0,
\end{equation}
where $\sigma_{s}^{kl}$ and $\sigma_{u}^{kl}$ are partial derivatives of $\sigma^{kl}$ with respect to $x^s$ and $u$ respectively. Since
\begin{equation}
\sigma_{\alpha}^{kl}=-(\sigma_{ij})_\alpha\sigma^{ik}\sigma^{jl},
\end{equation}
(\ref{0 order term}) becomes
\begin{equation}
\label{yellow}
\xi^s(\sigma_{kl})_{s}+\eta(\sigma_{kl})_{u}+\xi^{s}_{l}\sigma_{ks}+\xi^{s}_k\sigma_{sl}-\tau_t\sigma_{kl}=0
\end{equation}
All terms of the form $u^{p}u_{kl}$ are collected as
\begin{equation}
(\eta^{k}+\xi^{k}_{u})u^su_{sk}
\end{equation}
So
\begin{equation}
\eta^{k}+\xi^{k}_{u}=0.
\end{equation}
Lowering the indices, we get for each $k$
\begin{equation}
\label{orange}
\eta_{k}+\xi^{s}_{u}\sigma_{sk}=0.
\end{equation}

All terms of the form $u^iu^ju_{kl}$ can be divided into two parts: $u^ku^lu_{kl}$ and $u^iu^ju_{kl}$ for all $(i,j)\neq(k,l)$. The term of $u^ku^lu_{kl}$ is
\begin{equation}
(2\eta_u-\tau_t)u^ku^lu_{kl}
\end{equation}
If we can check that $u^ku^lu_{kl}$ and $u^iu^ju_{kl}$ for all $(i,j)\neq(k,l)$ are independent, then we get
\begin{equation}
\label{pink}
2\eta_u-\tau_t=0.
\end{equation}
We suppose otherwise there exist constants $C_0$ and $C_{pq}^{ij}$, which are not all zero, such that
\begin{equation}
\label{contra assump}
C_0u^ku^lu_{kl}+\sum_{(p,q)\neq(i,j)}C_{pq}^{ij}u^pu^qu_{ij}=0
\end{equation}
actually we can further assume that $C_{pq}^{ij}$ are symmetric both in $(i,j)$ and in $(p,q)$, since $u^pu^qu_{ij}$ have the same symmetries. Differentiating (\ref{contra assump}) with respect to $u_{ij}$ we get
\begin{equation}
C_0u^iu^j+\sum_{(p,q)\neq(i,j)}C_{pq}^{ij}u^pu^q=0
\end{equation}
Differentiating it further with respect to $u_{i}$ and by the symmetry of $C_{pq}^{ij}$ we get for $q\neq j$
\begin{equation}
C_0u^j+2\sum_{q\neq j}C_{iq}^{ij}u^iu^q=0
\end{equation}
Since $u_j$ and $u_q$ for $q\neq j$ are independent, we have $C_0=0$ and $C_{iq}^{ij}=0$, then (\ref{contra assump}) becomes for $p\neq i,j$ and $q\neq i,j$

\begin{equation}
\sum_{p\neq i,j\ q\neq i,j}C_{pq}^{ij}u^pu^qu_{ij}=0.
\end{equation}
Similarly we can get $C_{pq}^{ij}=0$, which is a contradiction. \par

In summary, we obtain the following system of equations from (\ref{yellow}), (\ref{orange}) and (\ref{pink})
\begin{equation}
\label{determ equ highdim}
\begin{cases}
&\xi^s(\sigma_{kl})_{s}+\eta(\sigma_{kl})_{u}+\xi^{s}_{l}\sigma_{ks}+\xi^{s}_k\sigma_{sl}=\tau_t\sigma_{kl}\\
&\eta_{k}+\xi^{s}_{u}\sigma_{sk}=0\\
&2\eta_u=\tau_t.
\end{cases}
\end{equation}
It can be checked that the left hand side of equations of (\ref{determ equ highdim}) are exactly the components of $L_Xg$, so they are tensor equations which is independent of the choice of coordinates of $M$, thus we can obtain that $\tau_{tt}=0$ by differentiating the first equation of (\ref{determ equ highdim}) with respect to $t$ and choosing a normal coordinate system $(x,u)$ around a point of $M$. Therefore, for some constant $\lambda$ the system (\ref{determ equ highdim}) can be written as
\begin{equation}
\label{conf kill equa highdim}
\begin{cases}
&\xi^s(\sigma_{kl})_{s}+\eta(\sigma_{kl})_{u}+\xi^{s}_{l}\sigma_{ks}+\xi^{s}_k\sigma_{sl}=2\lambda\sigma_{kl}\\
&\eta_{k}+\xi^{s}_{u}\sigma_{sk}=0\\
&\eta_u=\lambda.
\end{cases}
\end{equation}
It can be checked that (\ref{conf kill equa highdim}) is exactly $L_Xg=2\lambda g$ satisfied by the homothetically conformal vector field.\par
It remains to check the non-degeneracy condition. \par

\textbf{Step 3.} We first check the local solvability. By the same argument, we can assume without loss of generality that the initial data for the local solvability is
\begin{equation}
(x_0, 0, u_0, (u_t)_0, (u_i)_0, (u_{ij})_{0}, (u_{it})_{0}, (u_{tt})_{0}),
\end{equation}
for all $1\leqslant i, j \leqslant n$, where $x_0=(x_0^1,...,x_0^n)$, subject to the condition:
\begin{equation}
\begin{cases}
\Psi^1(x_0,t_0,u_0^{(2)})=0\\
\Psi^{kl}(x_0,t_0,u_0^{(2)})=0.
\end{cases}
\end{equation}
for all $k$ and $l$.

We are to construct a solution to the mean curvature flow $u=f(x,t)$ such that $u^{(2)}_{0}=pr^{(2)}f(x_0,0)$. We choose local normal coordinate system $\{x^1,...,x^n\}$ centered at the point $x_0$ of $M$, and define the initial hypersurface by
\begin{equation}
\begin{aligned}
f_0(x):=f(x,0)=&u_0+(u_i)_0(x^i-x_0^i)+\frac{1}{2}(u_{ij})_{0}(x^i-x_0^i)(x^j-x_0^j)\\
       &+\frac{1}{4!}C(x^1-x_0^1)^4,
\end{aligned}
\end{equation}
where $C$ is a constant to be determined. And $u=f(x,t)$ satisfies
\begin{equation}
\begin{cases}
\Psi^1(x,t,f^{(2)})=0\\
\Psi^{kl}(x,t,f^{(2)})=0\\
f(x,0)=f_0(x).
\end{cases}
\end{equation}

Using $\Psi^1(x_0,t_0,u_0^{(2)})=0$, and noting that $\Gamma^\alpha_{\beta\gamma}$ and $\sigma_{ij}$ depend only on $(x,u)$, we can obtain for all $i$ and $j$,
\begin{equation}
\begin{gathered}
f(x_0,0)=u_0,
\quad
f_i(x_0,0)=(u_i)_0,
\quad
f_{ij}(x_0,0)=(u_{ij})_0,
\quad
f_t(x_0,0)=(u_t)_0.
\end{gathered}
\end{equation}
By $\Psi^2(x_0,u_0^{(2)})=0$ and $\Psi^{kl}(x,f^{(2)})=0$, we see that for all $k$
\begin{equation}
f_{kt}(x_0,0)=(u_{kt})_0
\end{equation}

Since
\begin{equation}
\label{expr utt high dim}
u_{tt}=\frac{\partial L}{\partial t}H+L\frac{\partial H}{\partial t},
\end{equation}
and
\begin{equation}
\label{exp L}
\frac{\partial L^2}{\partial t}=2u_{kt}u_k,
\end{equation}
according to Lemma 10.7 of \cite{Zhu2002}, we have
\begin{equation}
\frac{\partial H}{\partial t}=\Delta H+H(|A|^2+Ric_N(\nu,\nu))
\end{equation}
Thus there exists a function $\alpha(x, u^{(3)})$, such that at the point $(x_0,u_0)$,
\begin{equation}
\label{exp Ht}
\begin{aligned}
\frac{\partial H}{\partial t}(x_0,u_0)=&\frac{1}{L_0}g_0^{ij}(f_0)_{ijkk}+\alpha(x_0,f_0^{(3)})\\
 =&\frac{C}{L_0}+\alpha(x_0,f_0^{(3)})
\end{aligned}
\end{equation}
where $L_0=L(x_0,u_0)$ and $g_0^{ij}=g^{ij}(x_0,u_0)$. Inserting (\ref{exp Ht}) and (\ref{exp L}) into (\ref{expr utt high dim}) and evaluating at the point $(x_0,u_0)$, since $L_0$ is positive, we can always choose proper $C$ such that $f_{tt}(x_0,u_0)=(u_{tt})_{0}$.\par

\textbf{Step 4.} Now we check the condition of maximal rank. If there is a point $(x_0,t_0,u_0^{(2)})\in \mathcal{S}:=\{(x,t,u^{(2)}):\Psi^1(x,t,u^{(2)})=0, \Psi^{kl}(x,t,u^{(2)})=0$, $1\leqslant k,l\leqslant n\}$ such that $J_{\Psi}$ is of maximal rank, then there is a neighbourhood of $(x_0,t_0,u_0^{(2)})$ such that $J_{\Psi}$ also have maximal rank due to smoothness of $\Psi$, so the condition is fulfilled in that neighbourhood where we can apply Theorem \ref{equi cond}. Thus, we only have to consider the case when there is no such a point, in another word, the determinants of all $(n^2+1)$ order minor matrices of $J_{\Psi}$ are zero for every point of $\mathcal{S}$.\par

We define
\begin{equation}
\Lambda_{(i,j)}^{(k,l)}:=\frac{1}{2L^2}\frac{\partial \Psi^{kl}}{\partial u_{ij}}
\end{equation}
and
\begin{equation}
\E_{i}^{(k,l)}:=\frac{1}{2L^2}\frac{\partial \Psi^{kl}}{\partial u_{it}}.
\end{equation}
The index $(k,l)$ and $(i,j)$ denote respectively the row and the column of the matrix $\Lambda_{(i,j)}^{(k,l)}$ with the lexicographical order. Then $\Lambda_{(i,j)}^{(k,l)}$ is a $n^2\times n^2$ matrix and $\E_{i}^{(k,l)}$ is $n^2\times n$ matrix.\par

We can calculate the elements of $\Lambda_{(i,j)}^{(k,l)}$ and $\E_{i}^{(k,l)}$ as follows.
\begin{equation}
\label{Lambda matrix}
\Lambda_{(i,j)}^{(k,l)}=
\begin{cases}
H+g^{ij}h_{ij}  &(k,l)=(i,j)\\
g^{ij}h_{kl}    &(k,l)\neq(i,j),
\end{cases}
\end{equation}
where there is no summation in $g^{ij}h_{ij}$.
\begin{equation}
\label{E matrix}
\E_{i}^{(k,l)}=
\begin{cases}
u_i    & k=l=i\\
\frac{u_l}{2}     &k\neq l, k=i \\
\frac{u_k}{2}     &k\neq l, l=i \\
0       & else.
\end{cases}
\end{equation}

For any $t\in I$, we can choose a local normal coordinates $\{x_1,...,x_n\}$ of some open subset $\mathcal{O}\subset M_t$ around any given point $p$, such that the metric $g(t)$ and the second fundamental form $h(t)$ of $M_t$ can be simultaneously diagonalized at $p$. Now we replace the column of $\Lambda_{(i,i)}^{(k,l)}$ with the corresponding the column of $\E_{i}^{(k,l)}$ for all $1\leqslant\ i \leqslant n$, and denote the new minor by $\tilde\Lambda_{(i,j)}^{(k,l)}$.\par

\begin{lemma}
\label{determinant lemma}
Each diagonal element of $\tilde\Lambda_{(i,j)}^{(k,l)}$ is the only non-zero element of the corresponding row or column.
\end{lemma}
\noindent
Proof. By (\ref{Lambda matrix}), (\ref{E matrix}) and the definition of $\tilde\Lambda_{(i,j)}^{(k,l)}$, we have
\begin{equation}
\label{(k,k)row}
\tilde\Lambda_{(i,j)}^{(k,k)}=
\begin{cases}
E_i^{(k,k)}=0,  & i=j, k\neq i\\
u_i,      & i=j, k=i\\
\Lambda_{(i,j)}^{(k,k)}=0, & i\neq j
\end{cases}
\end{equation}
and for $i\neq j$,
\begin{equation}
\label{i!=j column}
\tilde\Lambda_{(i,j)}^{(k,l)}=\Lambda_{(i,j)}^{(k,l)}
\begin{cases}
H,   & (k,l)=(i,j)\\
0,    & (k,l)\neq (i,j)
\end{cases}
\end{equation}
We note that the diagonal elements of $\tilde\Lambda$ are $\tilde\Lambda_{(i,j)}^{(k,l)}$ with $(k,l)=(i,j)$. By (\ref{(k,k)row}), the only non-zero element of $(k,k)$ row is the $\{(k,k),(k,k)\}$-th element which is $u_i$. By (\ref{i!=j column}), the only non-zero element of $(i,j)$ column for all $i\neq j$ is the $\{(k,l),(i,j)\}$-th element with $(k,l)=(i,j)$ which is $H$. Thus we complete the proof.\par
\hfill\qed

By Lemma \ref{determinant lemma} and Laplace expansion, the determinant of $\tilde\Lambda_{(i,j)}^{(k,l)}$ is the product of the diagonal elements
\begin{equation}
\det\{\tilde\Lambda_{(i,j)}^{(k,l)}\}=H^{n^2-n}\prod_{i=1}^{n}u_i.
\end{equation}
We note that
\begin{equation}
\frac{\partial \Psi^1}{\partial u_t}\neq 0
\end{equation}
and for all $k$ and $l$.
\begin{equation}
\quad \frac{\partial \Psi^{kl}}{\partial u_t}=0,
\end{equation}
By the assumption that each minor of order $n^2+1$ is degenerate on $\mathcal{O}$, we obtain
\begin{equation}
\det\{\tilde\Lambda_{(i,j)}^{(k,l)}\}=0.
\end{equation}
which implies $H=0$ or $u_i= 0$ on $\mathcal{O}$ for some $i$ at each time $t\in I$. Since the minimal hypersurface $H=0$ is a stationary solution to the mean curvature flow which admits only diffeomorphisms, it is a trivial case for our purpose. Thus we assume $H\neq0$ and $u_1=0$ without loss of generality.\par

Taking $k=l=1$ in the evolution equation of the metric, we get at the point $p$
\begin{equation}
0=2u_ku_{1t}=-2Hh_{11}=2H\mu_1,
\end{equation}
where $\mu_1$ is a principal curvature in the direction of $e_1$. We have $\mu_1=0$ due to $H\neq0$, and it is independent on our choice of the coordinates. Thus the $\{(1,1),(1,1)\}$-th element of the matrix $\Lambda_{(i,j)}^{(k,l)}$ is $H$, now we replace the column of $\Lambda_{(i,i)}^{(k,l)}$ with the column of $\E_{i}^{(k,l)}$ for $2\leqslant i\leqslant n$, and the determinant of the resulting matrix, which is also denoted by $\tilde\Lambda_{(i,j)}^{(k,l)}$, should also be zero by the degenerate assumption, that is
\begin{equation}
\det\{\tilde\Lambda_{(i,j)}^{(k,l)}\}=H^{n^2-n+1}\prod_{i=2}^{n} u_i.
\end{equation}
By the same arguments, we obtain $u_2=0$. Repeating the process, we obtain that $u_i=0$ for $1\leqslant i\leqslant n$, and the determinant of the original matrix $\Lambda_{(i,j)}^{(k,l)}$ becomes
\begin{equation}
\det\{\Lambda_{(i,j)}^{(k,l)}\}=H^{n^2}.
\end{equation}
Therefore, we obtain that when the solution is not a stationary minimal hypersurface, the system
\begin{equation}
\begin{cases}
\Psi^1(x,t,u^{(2)})=0\\
\Psi^{kl}(x,t,u^{(2)})=0
\end{cases}
\end{equation}
is of maximal rank.\par
\hfill\qed

\section{Mean curvature flow solitons}

In this section, we derive the characterizing equation of the mean curvature flow solitons and give examples of the homothetic solitons in non-Euclidean surfaces.\par

We first give the definition of the mean curvature flow solitons
\begin{defn}
A smooth solution $F :M^{n}\times I \rightarrow (N^{n+1},\bar g)$ to the mean curvature flow is called a soliton if there exists a one-parameter subgroup $\omega_t$, $t\in I$ of the symmetry group of the mean curvature flow, such that $F(x,t)=\omega_t\cdot F_{0}(x)$, where $x\in M$ and $F_{0}(x)=F(x,0)$.
\end{defn}

Let $F(x,t)$ be the mean curvature flow soliton, such that
\begin{equation}
\begin{cases}
\big(\frac{\partial F}{\partial t}\big)^\bot=H\nu\\
F(x,0)=F_0(x)
\end{cases}
\end{equation}
and we assume that $\omega_t$ is a local one-parameter homothetic transformations of $N$ satisfying $\omega_{t}^*(\bar g)=c^2(t)\bar g$ for some positive function $c(t)$, and
\begin{equation}
\label{one-parameter equ}
\begin{cases}
\frac{d\omega_t(x)}{dt}=X(\omega_t(x))\\
\omega_0(x)=x.
\end{cases}
\end{equation}
We call $X$ the homothetic vector field corresponding to the one-parameter homothetic transformations, and we have the following relations
\begin{equation}
L_X\omega^*_{t}(\bar g)=\lim_{\varepsilon \to 0}\frac{\omega^*_{t+\varepsilon}(\bar g)-\omega^*_{t}(\bar g)}{\varepsilon}=2c'(t)c(t)\bar g.
\end{equation}

By the definition of solitons,  we have
\begin{equation}
\bigg(\frac{\partial F}{\partial t}\bigg)^\bot(\omega(x,t),t)=X^\bot(F(x,t))=H(F(x,t))\nu(F(x,t))
\end{equation}
\noindent
Applying tangential mappings $(\omega_t^{-1})_*=(\omega_{-t})_*$ to the above equation, we obtain
\begin{equation}
(\omega_{-t})_*(X^\bot(F(x,t))=H(F(x,t))(\omega_{-t})_*(\nu(F(x,t))).
\end{equation}
Since $(\omega_t)_*$ is an isomorphism between $T_{F_0(x)}M$ and $T_{F(x,t)}M$, in addition, the conformal mapping preserves the orthogonality, we have
\begin{equation}
(\omega_{-t})_*(X^\bot(F(x,t))=\big[(\omega_{-t})_*(X(F(x,t))\big]^\bot
\end{equation}
and we can also fix an orientation of the hypersurface such that
\begin{equation}
(\omega_{-t})_*(\nu(F(x,t)))=\frac{1}{c(t)}\nu(x,0).
\end{equation}
The corresponding mean curvatures can be related by
\begin{equation}
H(F(x,t))=\frac{1}{c(t)}H(F_0(x)).
\end{equation}
Therefore the characterizing equation of the mean curvature flow solitons is
\begin{equation}
\label{character equ general}
\big[(\omega_{-t})_*(X(F(x,t))\big]^\bot=\frac{1}{c^2(t)}H(x,0)\nu(x,0).
\end{equation}
Particularly, when $t=0$, we have
\begin{equation}
\label{charecter equ}
X^\bot=H\nu
\end{equation}

As an example we recover the classical homothetic solutions in Euclidean spaces.\par
\begin{example}
The homothetic transformation of $\mathbb{R}^{n+1}$ is defined by
\begin{equation}
\omega_t(x)=c(t)x,
\end{equation}
for $x=(x^1,...,x^{n+1})\in\mathbb{R}^{n+1}$, and a homothetic soliton is $F(x,t)=c(t)F_0(x)$, for some given initial hypersurface $F_0:M \rightarrow \mathbb{R}^{n+1}$. The corresponding vector field $X$ of $\omega_t(x)$ is then
\begin{equation}
X(F(x,t))=c'(t)x(F(x,t))
\end{equation}

Since $\omega_t^{-1}=\frac{1}{c(t)}x$, we have
\begin{equation}
(\omega_t^{-1})_*=\frac{1}{c(t)}Id
\end{equation}
then
\begin{equation}
(\omega_t^{-1})_*(X(F(x,t)))=\frac{c'(t)}{c(t)}x(\omega_t(F(x,t)))=\frac{c'(t)}{c(t)}F_0(x).
\end{equation}
Thus by (\ref{character equ general}), we get
\begin{equation}
\big[(\omega_t^{-1})_*(X(F(x,t)))\big]^\bot=\frac{1}{c^2(t)}H(x,0)\nu(x,0)=\frac{c'(t)}{c(t)}F_0^\bot(x),
\end{equation}
that is,
\begin{equation}
H(x,0)\nu(x,0)=c(t)c'(t)F_0^\bot(x).
\end{equation}
This is exactly the characterizing equation of the homothetic solutions in Euclidean spaces see \cite{Eck2004} for further accounts.
\end{example}

We next construct examples of non-Euclidean homothetic solitons generated by a special one-parameter homothetic transformations $\omega_t$ defined by (\ref{one-parameter equ}), where we further assume that $L_X\bar g=2\lambda\bar g$ for some non-zero constant $\lambda$. From the following lemmas, we can see that there is an obstruction to the existence of such homothetic transformations.

\begin{lemma}
\label{kob lemma}
(Kobayashi \cite{Kob1963}, pp. 242, Lemma 2)
If $M$ is a complete Riemannian manifold which is not locally Euclidean, then any homothetic transformation of $M$ is an isometry.
\end{lemma}

\begin{lemma}
\label{homo lemma}
If a complete Riemannian manifold $M$ admits one-parameter homothetic transformations defined above, then $M$ is isometric to a Euclidean space.
\end{lemma}
\noindent
Proof. By Lemma \ref{kob lemma}, $M$ must be locally Euclidean and thus a flat manifold. Since $M$ is a complete flat manifold, its universal Riemannian covering space is a Euclidean space (cf. \cite{Kob1963}). Let $\tilde V$ be the homothetic vector field induced by the the family of transformations and $V$ be its horizontal lift. Since Riemannian covering is local isometry, the vector field $V$ is a homothetic vector field on $(\mathbb{R}^{n},g_0)$, where $g_0=\sum_{i=1}^{n}(dx^i)^2$. Suppose $V=\sum_{i=1}^{n}V^i\frac{\partial}{\partial x^i}$, such that $L_Vg_0=2\lambda g_0$, where $\lambda$ is non-zero constant, then we have

\begin{equation}
g_0(D_{\frac{\partial}{\partial x^i}}{V},\frac{\partial}{\partial x^j})+g_0( D_{\frac{\partial}{\partial x^j}}V,\frac{\partial}{\partial x^i})=2\lambda \delta_{ij}.
\end{equation}
which can be further reduced to
\begin{equation}
V_i^j+V_j^i=2\lambda\delta_{ij},
\end{equation}
where $V_i^j=\frac{\partial V^j}{\partial x^i}$. We see that $V_i^j$ satisfies
\begin{equation}
\begin{cases}
V_i^i=\lambda,  & 1\leqslant i \leqslant n\\
V_i^j=-V_j^i,   & 1\leqslant i\neq j \leqslant n.
\end{cases}
\end{equation}
Thus for each $i$, there exists a function $\alpha^i$ satisfying
\begin{equation}
\begin{cases}
\frac{\partial \alpha^i}{\partial x^i}=0,        & 1\leqslant i \leqslant n\\
\frac{\partial \alpha^i}{\partial x^j}=-\frac{\partial \alpha^j}{\partial x^i},  & 1\leqslant i\neq j \leqslant n.
\end{cases}
\end{equation}
such that
\begin{equation}
V^i(x)=\lambda x^i+\alpha^i(x^1,...\widehat{x^i}...,x^n)
\end{equation}
where $\widehat{x^i}$ represents omitting the variable. \par

If we define
\begin{equation}
W:=\sum_{i=1}^{n}\alpha^i(x^1,...\widehat{x^i}...,x^n)\frac{\partial}{\partial x^i},
\end{equation}
then it is straightforward to check that $L_Wg_0=0$, so $W$ is a Killing vector field. Since a Killing vector field on $\mathbb{R}^n$ can be expressed by the linear combination of vector fields corresponding to translations and rotations, we have
\begin{equation}
W=\sum_{i=1}^{n}T^i\frac{\partial}{\partial x^i}+\sum_{p\neq q}^{n}R_{pq}\{-x^p\frac{\partial}{\partial x^q}+x^q\frac{\partial}{\partial x^p}\},
\end{equation}
and we can choose $R_{pq}$ to be an skew symmetric matrix without loss of generality.\par

We assume the covering map is
\begin{equation}
\pi:\mathbb{R}^n\rightarrow \mathbb{R}^n/\Gamma
\end{equation}
where $\Gamma$ is a lattice, namely
\begin{equation}
\Gamma=Span_{\mathbb{Z}}\{\omega_1,\ ...,\ \omega_s\}
\end{equation}
where $\omega_1$, ..., $\omega_s$, $1\leqslant s \leqslant n$ are independent vectors of $\mathbb{R}^{n}$. Since $\tilde V=(\pi)_*(V)$ is global vector field, $V$ must be
$\Gamma-$periodic, that is
\begin{equation}
\label{period}
V(x)=V(x+k\omega_i),
\end{equation}
for arbitrary $k\in\mathbb{Z}$ and $1\leqslant i \leqslant s$. For any non-zero vector $\omega\in\Gamma$, supposing its coordinates are $(a^1,...,a^n)$, we get by (\ref{period}),
\begin{equation}
\begin{aligned}
\sum_{i=1}^{n}\lambda x^i\frac{\partial}{\partial x^i}+W=&\sum_{i=1}^{n}\lambda(x^i+ka^i)\frac{\partial}{\partial x^i}+\sum_{i=1}^{n}T^i\frac{\partial}{\partial x^i}\\
       &+\sum_{p\neq q}^{n}R_{pq}\{-(x^p+ka^p)\frac{\partial}{\partial x^q}+(x^q+ka^q)\frac{\partial}{\partial x^p}\}.
\end{aligned}
\end{equation}
By the skew symmetry of $R:=R_{pq}$, we obtain
\begin{equation}
(\lambda-2R)\omega=0.
\end{equation}
since $\omega\neq0$, $\lambda$ is a eigenvalue of $2R$. Since the only real eigenvalue of an skew symmetric matrix is zero, it contradicts with the assumption.\par
\hfill\qed

We now give examples of homothetic solitons on surface patches.

\begin{example}
We consider a Riemannian surface patch $(\Sigma, g)$, where $\Sigma$ is a connected domain with isothermal coordinates $\{u,v\}$, such that
\begin{equation}
g=e^{2\rho(u,v)}(du^2+dv^2).
\end{equation}
Suppose also that
\begin{equation}
X=\phi(u,v)\frac{\partial}{\partial u}+\psi(u,v)\frac{\partial}{\partial v}
\end{equation}
is a conformal vector field on $\Sigma$ such that $L_Xg=2\lambda g$ for some constant $\lambda$, and we are particularly interested in the case when $\lambda\neq0$. By straightforward calculation, the condition $L_Xg=2\lambda g$ implies that $\phi$ and $\psi$ satisfies
\begin{equation}
\begin{cases}
\phi_u+\phi\rho_u+\psi\rho_v=\lambda\\
\phi_u-\psi_v=0\\
\phi_v+\psi_u=0.
\end{cases}
\end{equation}

Suppose $r:I\rightarrow \Sigma$ is a smooth curve, where $I$ is an open interval, and we choose $s$ to be the arc-length parameter. Then the tangent field of $r(s)=(u(s),v(s))$ is
\begin{equation}
r_*(\frac{\partial}{\partial s})=u'\frac{\partial}{\partial u}+v'\frac{\partial}{\partial v}
\end{equation}
where $u'$ and $v'$ are partial derivatives with respect to $s$, satisfying
\begin{equation}
\big|r_*(\frac{\partial}{\partial s})\big|=1,
\end{equation}
which implies
\begin{equation}
\label{length 1 relation}
(u')^2+(v')^2=e^{-2\rho}.
\end{equation}
The left-ward pointing unit normal $\nu$ is then
\begin{equation}
\nu=-v'\frac{\partial}{\partial u}+u'\frac{\partial}{\partial v}.
\end{equation}
By Gauss formula we obtain the geodesic curvature $k_g$ of $r$
\begin{equation}
k_g=e^{2\rho}(-v'f_1+u'f_2),
\end{equation}
where
\begin{equation}
\begin{cases}
f_1:=u''+(u')^2\rho_u-(v')^2\rho_u+2u'v'\rho_v\\
f_2:=v''-(u')^2\rho_v+(v')^2\rho_v+2u'v'\rho_u
\end{cases}
\end{equation}
By the characterizing equation of the mean curvature flow soliton (\ref{charecter equ}), we obtain
\begin{equation}
(\phi-f_1)v'+(f_2-\psi)u'=0.
\end{equation}
Therefore in order to determine the mean curvature flow solitons on the surface $(\Sigma,g)$, we have to solve the following system of equations
\begin{subequations}
\label{example equ}
\begin{empheq}[left=\empheqlbrace]{align}
&(\phi-f_1)v'+(f_2-\psi)u'=0\\
&\phi_u+\phi\rho_u+\psi\rho_v=\lambda\\
&\phi_u-\psi_v=0\\
&\phi_v+\psi_u=0.
\end{empheq}
\end{subequations}
Some special solutions to this system are easily obtained.\par
(I) Firstly, we observe that a special type of solutions to (\ref{example equ}c) and (\ref{example equ}d) are
\begin{equation}
\begin{cases}
\phi=au+bv+c_1\\
\psi=-bu+av+c_2.
\end{cases}
\end{equation}
where $a,b,c_1$ and $c_2$ are constants. Secondly, we observe that (\ref{example equ}b) is a first order linear partial differential equation with respect to $\rho$, thus we can solve it in the following cases.

(\romannumeral1) $a=b=0$.
The characteristic ODE system of (\ref{example equ}b) is
\begin{equation}
\begin{cases}
\frac{du}{dt}=c_1\\
\frac{dv}{dt}=c_2\\
\frac{d\rho}{dt}=\lambda.
\end{cases}
\end{equation}
and its first integrals are
\begin{equation}
I_1:=c_2u-c_1v \quad and \quad I_2:=\lambda u-c_1\rho
\end{equation}
or
\begin{equation}
I_1:=c_2u-c_1v \quad and \quad I_2:=\lambda v-c_2\rho
\end{equation}
the corresponding solutions of $\rho$ are
\begin{equation}
\rho=\frac{\lambda}{c_1}u+Q(c_2u-c_1v)
\end{equation}
or
\begin{equation}
\rho=\frac{\lambda}{c_2}v+Q(c_2u-c_1v)
\end{equation}
where $Q$ is any smooth function.\par
Now the Gauss curvature $K$ of $\Sigma$ can be computed by
\begin{equation}
K=-e^{-2\rho}(\rho_{uu}+\rho_{vv}).
\end{equation}
which is generally not zero, we thus obtain non-Euclidean homothetic solitons when $\lambda\neq 0$. \par

Differentiating (\ref{length 1 relation}) and combining (\ref{example equ}a), we obtain the ODE system satisfied by the initial curve $r$
\begin{equation}
\label{curve ODEs}
\frac{d}{ds}
\left[
\begin{array}{c}
   u\\
   v\\
   w\\
   z.
\end{array}
\right]
=
\left[
\begin{array}{c}
   w\\
   z\\
   h_1(u,v,w,z)\\
   h_2(u,v,w,z)
\end{array}
\right]
\end{equation}
where
\begin{equation}
\begin{aligned}
h_1(u,v,w,z):=&-\frac{1}{w}\{zz'+e^{-2\rho}(w\rho_u+z\rho_v)\}\\
h_2(u,v,w,z):=&\frac{z}{w}(w'+w^2\rho_u-z^2\rho_u+2wz\rho_v-c_1)\\
              &+w^2\rho_v-z^2\rho_v-2wz\rho_u+c_2.
\end{aligned}
\end{equation}

(\romannumeral2) $b=0$ and $a\neq0$. Now (\ref{example equ}b) becomes
\begin{equation}
(au+c_1)\rho_u+(av+c_2)\rho_v=\lambda-a.
\end{equation}
Its characteristic ODE system is
\begin{equation}
\begin{cases}
\frac{du}{dt}=au+c_1\\
\frac{dv}{dt}=av+c_2\\
\frac{d\rho}{dt}=\lambda-a.
\end{cases}
\end{equation}
By
\begin{equation}
\frac{d\rho}{\lambda-a}=\frac{du}{au+c_1},
\end{equation}
we obtain a special solution
\begin{equation}
\rho=\frac{\lambda-a}{a}\ln|au+c_1|.
\end{equation}
Multiplying
\begin{equation}
\frac{du}{au+c_1}-\frac{dv}{av+c_2}=0
\end{equation}
with an integrating factor
\begin{equation}
\mu=\frac{a^2|au+c_1|^a}{|av+c_2|^a},
\end{equation}
we can get
\begin{equation}
d(|au+c_1|^a|av+c_2|^{-a})=0,
\end{equation}
and it implies that $I_1=|au+c_1|^a|av+c_2|^{-a}$ is a first integral. Therefore the solution can be written as
\begin{equation}
\rho=\frac{\lambda-a}{a}\ln|au+c_1|+Q(|au+c_1|^a|av+c_2|^{-a}),
\end{equation}
for $u\neq-\frac{c_2}{a}$, $v\neq-\frac{c_1}{a}$ and an arbitrary smooth function $Q$.\par
On a connected sub-domain
\begin{equation}
D\subset \{(u,v)\in \Sigma: u\neq-\frac{c_2}{a},v\neq-\frac{c_1}{a}\},
\end{equation}
a soliton can be solved by (\ref{curve ODEs}) defined on $D$.\par

(\romannumeral3) $a=0$ and $b\neq 0$. Now (\ref{example equ}b) is
\begin{equation}
(bv+c_1)\rho_u+(-bu+c_2)\rho_v=\lambda.
\end{equation}
and its characteristic ODE system is
\begin{equation}
\begin{cases}
\frac{du}{dt}=bv+c_1\\
\frac{dv}{dt}=-bu+c_2\\
\frac{d\rho}{dt}=\lambda.
\end{cases}
\end{equation}
By
\begin{equation}
\frac{du}{bv+c_1}=\frac{d\rho}{\lambda},
\end{equation}
we obtain a special solution
\begin{equation}
-\frac{\lambda}{b}\arctan\bigg(\frac{bv+c_1}{bu-c_2}\bigg),
\end{equation}
where we assume that $u\neq\frac{c_2}{b}$ and $v\neq-\frac{c_1}{b}$.
By
\begin{equation}
\frac{du}{bv+c_1}=\frac{dv}{-bu+c_2}
\end{equation}
we obtain a first integral
\begin{equation}
I_1:=(bv+c_1)^2+(bu-c_2)^2.
\end{equation}
Thus the general solution can be written as
\begin{equation}
\rho=-\frac{\lambda}{b} \arctan\bigg(\frac{bv+c_1}{bu-c_2}\bigg)+Q[(bv+c_1)^2+(bu-c_2)^2],
\end{equation}
for $u\neq\frac{c_2}{b}$, $v\neq-\frac{c_1}{b}$ and an arbitrary smooth function $Q$. \par

Similar to the example in case (II), a soliton solution can be constructed on some connected sub-domain of $\Sigma$.\par

(II) A homogeneous polynomial solution to (\ref{example equ}c) and (\ref{example equ}d) is
\begin{equation}
\begin{cases}
\phi=u^2-v^2\\
\psi=2uv
\end{cases}
\end{equation}
and the corresponding characteristic ODE system is
\begin{equation}
\begin{cases}
\frac{du}{dt}=u^2-v^2\\
\frac{dv}{dt}=2uv\\
\frac{d\rho}{dt}=\lambda-2u.
\end{cases}
\end{equation}
Since
\begin{equation}
2uvdu+(v^2-u^2)dv=0,
\end{equation}
we introduce the following coordinates transformation
\begin{equation}
\begin{cases}
u=r\sinh\theta\\
v=r\cosh\theta,
\end{cases}
\end{equation}
and obtain
\begin{equation}
\frac{dr}{r}=-\frac{z^3+3z^2-3z-1}{2z(z+1)(z^2+1)}dz,
\end{equation}
where
\begin{equation}
z:=e^{2\theta}=\frac{u+v}{v-u}.
\end{equation}
From this equation $r$ can be expressed by $z$ as follows
\begin{equation}
r=C\frac{\sqrt z (z+1)}{z^2+1},
\end{equation}
for any positive constant $C$. Inserting $r$ into
\begin{equation}
\frac{d\rho}{dz}=\frac{\lambda-2u}{r^2}\frac{du}{dz},
\end{equation}
we can solve $\rho(z)$
\begin{equation}
\rho(z)=\frac{4\lambda}{c(z+1)}+4\ln\frac{(z+1)^2}{z^2+1}+D,
\end{equation}
where $D$ is an arbitrary constant.
\end{example}

\section{Proof of Theorem 2}
As an application of Theorem 1, we consider the affine solutions to the mean curvature flow in the Euclidean space $\mathbb{R}^{n}$.\par

\begin{defn}
A solution $F(x,t)$ to the mean curvature flow is called an affine solution, if there exists a one-parameter family of affine transformations $A(t)$ of $\mathbb{R}^{n}$,  such that $M_t=A(t)(M_0)$.
\end{defn}
Proof of Theorem 2.
In a Cartesian coordinate system $\{y^1,...,y^{n}\}$ of $\mathbb{R}^{n}$, the affine transformations can be represented as follows
\begin{equation}
F(x,t)=R(t)F_0(x)+T(t),
\end{equation}
where $T(t)$, $F(x,t)$ and $F_0(x)$ are regarded as column vectors and $R(t)$ is $n\times n$ matrix. Here we assume that the determinant of $R(t)$ is positive for any $t\in I$. By the initial condition $A(0)=Id$, we see that $R(0)$ is a unit matrix $I_n$, and $T(0)=0$.\par

According to QR decomposition, a matrix can always be decomposed into the product of an orthogonal matrix and an upper triangular matrix. Thus we assume there exist an orthogonal matrix $\tilde U(t)$ and an upper triangular matrix $\tilde V(t)$, such that
\begin{equation}
R(t)=\tilde U(t)\tilde V(t).
\end{equation}
Since det$\{R(t)\}>0$, the diagonal elements of $\tilde V(t)$ are non-zero. Then we can further have the following decomposition
\begin{equation}
R(t)=
\left[
\begin{array}{ccc}
u_{11}(t)& \cdots  & u_{1n}(t) \\
\vdots   & \ddots  & \vdots    \\
u_{n1}(t) & \cdots & u_{nn}(t)
\end{array}
\right]
\left[
\begin{array}{ccc}
s(t)     &        &0   \\
         & \ddots &   \\
0        &        & s(t)
\end{array}
\right]
\left[
\begin{array}{ccc}
1        & \cdots & v_{1n}(t) \\
         & \ddots &\vdots     \\
0        &        & v_{nn}(t)
\end{array}
\right]
\end{equation}
and we denote it by
\begin{equation}
R(t)=U(t)S(t)V(t),
\end{equation}
where $U(t)$ is a special orthogonal matrix, $S(t)$ is a scalar matrix and $V(t)$ is an upper triangular matrix, with $U(0)=S(0)=V(0)=I_n$.

Now we define a one-parameter family of transformations acting on $F(x,t)$ with $t\in[t_0,t_2]\subset I$, such that they are still solutions to the mean curvature flow after the action. Let $\varepsilon\in[-\delta, \delta]$, where $\delta$ is small enough such that $t+\varepsilon\in I$, the action is defined by
\begin{equation}
A(\varepsilon)(F(x,t))=F(x,t+\varepsilon),
\end{equation}
where $t\in[t_0,t_1]$. In another word, $A(\varepsilon)$ is a one-parameter subgroup of the symmetry group. Since $A(t)$ is a one-parameter family of transformations, it satisfy $A(t+\varepsilon)=A(\varepsilon)A(t)$, thus
\begin{equation}
F(x,t+\varepsilon)=A(t+\varepsilon)(F_0(x))=A(\varepsilon)A(t)(F_0(x))=R(\varepsilon)[R(t)F_0(x)+T(t)]+T(\varepsilon).
\end{equation}
Differentiating $F(x,t+\varepsilon)$ with respect to $\varepsilon$ yields
\begin{equation}
\frac{\partial F(x,t+\varepsilon)}{\partial \varepsilon}\bigg|_{\varepsilon=0}=R'(0)F(x,t)+T'(0),
\end{equation}
and noting that
\begin{equation}
R'(0)=U'(0)+S'(0)+V'(0),
\end{equation}
we have
\begin{equation}
\begin{aligned}
X:=&\frac{\partial F(x,t)}{\partial t}=\frac{\partial F(x,t+\varepsilon)}{\partial \varepsilon}\bigg|_{\varepsilon=0}\\
=&U'(0)(F(x,t))+S'(0)(F(x,t))+V'(0)(F(x,t))+T'(0).
\end{aligned}
\end{equation}
Recall that $so(n)=\{A\in gl(n):A^t+A=0\}$, we see that $V'(t)$ is not in $so(n)$ unless $V(t)$ is a unit matrix. Also $V(t)$ is not scalar matrix unless $V(t)$ is a unit matrix. So $V(t)$ is generally neither a rotation nor a scaling and obviously not a translation. Since $X$ is an infinitesimal symmetry, by Theorem 1, $V'(t)$ must be tangent vector field on $M$, so $V(t)$ is a family of diffeomorphisms of $M$. If we consider the mean curvature flow with only normal motion, then it is only possible that $U(t)=I$. If we consider the general mean curvature flows, then by $X^{\bot}=H\nu$, they are exactly the self-similar solutions combining translation, rotation and scaling. And this proves Theorem 2.\par
\hfill\qed

\bibliographystyle{abbrv}

\begin{thebibliography}{10}


\bibitem{AbL1986}
U. ~Abresch and J. ~Langer.
\newblock The normalized curve shortening flow and homothetic solutions.
\newblock {\em J. Differ. Geom.}, 23, 175-196, 1986.

\bibitem{ALR2020}
J. ~Alias, J. de ~Lira, and M. ~Rigoli.
\newblock Mean curvature flow solitons in the presence of conformal vector fields.
\newblock {\em J. Geom. Anal.}, 30, 1466-1529, 2020.

\bibitem{Alt1991}
S. ~Altschuler.
\newblock Singularities of the curve shrinking flow for space curves.
\newblock {\em J. Differ. Geom.}, 34, 491-514, 1991.

\bibitem{Cho2002}
K. ~Chou and G. ~Li.
\newblock Optimal systems of group invariant solutions for the generalized curve shortening flow.
\newblock {\em Communications In Analysis And Geometry}, Volume 10, 241-274, 2002.

\bibitem{Eck2004}
K. ~Ecker.
\newblock Regularity theory for mean curvature flow.
\newblock Birkhaeuser, Boston-Basel-Berlin, 2004.

\bibitem{Fut2014}
A. ~Futaki, K. ~Hattori and H. ~Yamamoto.
\newblock Self-similar solutions to the mean curvature flows on Riemannian cone manifolds and special Lagrangians on toric Calabi-Yau cones.
\newblock {\em Osaka J. Math.}, 51, 1053-1079, 2014.

\bibitem{GaH1986}
M. ~Gage and R. ~Hamilton.
\newblock The heat equation shrinking convex plane curves.
\newblock {\em J. Differ. Geom.}, 23, 417-491, 1986.

\bibitem{Ger2006}
C. ~Gerhardt.
\newblock Curvature problems.
\newblock Series in geometry and topology, volume 39, International Press, 2006.

\bibitem{Hal2012}
H. ~Halldorsson.
\newblock Self-similar solutions to the curve shortening flow.
\newblock {\em Trans. Amer. Math. Soc.}, 364, 5285-5309, 2012.

\bibitem{Ham1995}
R. ~Hamilton.
\newblock Harnack estimate for the mean curvature flow.
\newblock {\em J. Differ. Geom.},  41, 215-226, 1995.

\bibitem{Hui1990}
G. ~Huisken.
\newblock Asymptotic behaviour for singularities of the mean curvature flow.
\newblock {\em J. Differ. Geom.}, 31, 285-299, 1990.

\bibitem{Hun2000}
N. ~Hungerbuehler and K. ~Smoczyk.
\newblock Soliton solutions for the mean curvature flow.
\newblock {\em Differential and Integral Equations}, Volume 13, 1321-1345, 2000.

\bibitem{Kob1963}
S. ~Kobayashi.
\newblock Foundations of differential geometry.
\newblock Vol 1, Interscience Publisher, 1963.

\bibitem{Olv1993}
P. ~Olver.
\newblock Applications of Lie groups to the differential equations.
\newblock Springer, 2nd Edition, 1993.

\bibitem{OST1994}
P. ~Olver, G. ~Sapiro and A. ~Tannenbaum.
\newblock Classification and uniqueness of invariant geometric flows.
\newblock {\em C. R. Acad. Sci. Paris Sér. I Math.}, 319(4), 339-344, 1994.

\bibitem{Smk2001}
K. ~Smoczyk.
\newblock A relation between mean curvature flow solitons and minimal submanifolds.
\newblock {\em Math. Nachr.}, 229, 175-186, 2001.

\bibitem{Zhu2002}
X. ~Zhu.
\newblock Lectures on mean curvature flows.
\newblock {\em International Press}, studies in advanced mathematics, vol. 32, 2002.




\end{thebibliography}

\end{document}